\documentclass[12pt,leqno,english]{article}
\usepackage{amsmath}
\usepackage{amsfonts}
\usepackage{amssymb}

\newcommand{\R}{\mathbb{R}}
\newcommand{\T}{\mathbb{T}}
\newcommand{\Z}{\mathbb{Z}}

\newcommand{\FF}{\mathcal{F}}

\newcommand{\ba}{\backslash}
\newcommand{\rg}{\rightarrow}

\begin{document}

\title
{
A tribute to the memory of Paul Tur\'an\\
\vskip2cm
Tur\'an's New Method and Compressive Sampling}

\author{Jean--Pierre Kahane}

\date{}

\maketitle

Tur\'an's New Method of Analysis of \cite{Turan} is a systematic study of trigonometric polynomials, from different points of view, with very different methods, and with a vast amount of applications. My purpose is to show a new application in the domain of signal's theory. I shall start with a theorem (theorem 11.6 in \cite{Turan}) and a question asked about, explain how to improve the theorem and discuss the question in relation with other results of \cite{Turan}: it will be parts~I and II of this article. Parts~III and IV will show its relation with a recent result of mine, a variation about a theorem of Cand\`es, Romberg and Tao of great interest in compressive sampling (theorem 1.3 in \cite{CaRoTao} and theorem~2.1 in~\cite{Cand}.

\section{The theorem and its improvement}

\textsc{Theorem} A (11.6 in \cite{Turan}).--- \textit{Given an integer $n\ge2$ and $0<\delta<1$, there exist real $x_1,x_2,\ldots, x_n$ such that}
$$
\Big| \sum_{j=1}^n e^{2\pi imx_j} \Big| < \delta n
\leqno(1)
$$ 
\textit{for all integers $m$ such that}
$$
1 \le m \le \frac{1}{2} 2^{n\delta^{2}/4}
\leqno(2)
$$

\vskip2mm

The proof relies on two ideas : 1) The distribution of the $x_j$ modulo 1 is the same as the distribution of the $mx_j$ modulo~1. 2) Considering the $x_j$ as independent random variables uniformly distributed on $\T = \R/\Z$, there exists an explicit formula for the moments of $|\sum\limits_{j=1}^n e^{2\pi ix_j}|$, namely
$$
E\Big(\Big| \sum_{j=1}^n e^{2\pi ix_j}\Big|^{2p}\Big) = \sum_{p_1+\cdots+p_k=2p}\ \frac{(2p)!}{p_1 !\cdots p_k !}\,.
$$
A preliminary lemma (5.9 in \cite{Turan}) gives an estimate of the right--hand member and a rather simple calculation leads to the result.

The improvement consists in using Laplace transforms. We are looking~for, taking into account that the $\cos (2 \pi x_j - \phi)$ and the  $\cos 2 \pi x_j$  have the same distribution, 
$$
P\Big(\Big| \sum_{j=1}^n e^{2\pi i x_j}\Big| > \delta n\Big)
$$
Let us choose an integer $\nu \ge 3$. Given any $z$ complex $\not=0$, there exists a $\varphi= \frac{2j\pi}{\nu}$ $(j=1,2,\ldots \nu)$ such that $Re\, ze^{-i\varphi}\ge a|z|$ with $a=\cos \frac{\pi}{\nu}$. Therefore
$$
\begin{array}{ll}
P\Big(\Big| \displaystyle \sum_{j=1}^n e^{2\pi ix_j}\Big| > \delta n\Big)
	&\displaystyle \le \nu \sup_\varphi P\Big(\sum_{j=1}^n \cos (2\pi x_j -\varphi)>\delta na\Big)\\
	\noalign{\vskip2mm}
	&\displaystyle =\nu P\Big(\sum_{j=1}^n \cos 2\pi x_j >\delta na\Big)\\
	\noalign{\vskip2mm}
	&\displaystyle <\nu E\Big(\exp u\Big(\sum_{j=1}^n \cos 2\pi x_j -\delta na\Big)\Big) \qquad (u>0)\\
	\noalign{\vskip2mm}
	&\displaystyle =\nu e^{-u \delta na}\prod_{j=1}^nE( e^{u\cos 2\pi x_j})\\
	\noalign{\vskip2mm}
	&= \nu e^{-u\delta na} (Ee^{u\cos 2\pi x}))^n
\end{array}
$$
where $x$ is a random variable uniformly distributed on $\T$. Now
$$
\begin{array}{ll}
E(e^{u \cos 2\pi x}) 
	&\displaystyle = \sum_{j=0}^\infty \frac{u^j}{j!} E\,((\cos 2\pi x)^j)\\
	\noalign{\vskip2mm}
	&=\displaystyle \sum_{k=0}^\infty \frac{u^{2k}}{(2k)!} \frac{(2k)!}{(k!)^2\,2^{2k}}\\
	\noalign{\vskip2mm}
	&\displaystyle\le\sum_{k=0}^\infty \frac{1}{k!} \Big(\frac{u^2}{4}\Big)^k = e^{+u^{2}/4}\\
\end{array}
$$
Therefore
$$
\begin{array}{ll}
P\Big(\Big| \displaystyle \sum_{j=1}^n e^{2\pi ix_j}\Big| > \delta \Big)
	& \displaystyle <\nu \inf_{u>0} e^{-u\delta na+n u^{2}/4}\\
\noalign{\vskip2mm}
	&=\nu e^{-n\delta^2a^2} \qquad \Big(a^2 =\cos^2 \dfrac{\pi}{\nu}\Big)	
\end{array}
$$

The best choice of $\nu$ minimizes $\log \nu -n \delta^2 \cos^2 \frac{\pi}{\nu}$. Since
$$
\frac{d}{dv} (\log v - n \delta^2 \cos^2 \frac{\pi}{v}) = \frac{1}{v}\Big(1-n \delta^2 \frac{\pi}{v} \sin \frac{2\pi}{v}\Big)
$$
we are lead to choose
$$
\nu= [\delta\pi \sqrt{2n}]
$$

Finally we obtain :

\vskip2mm

\textsc{Theorem} 1.--- \textit{Considering the $x_j$ $(j=1,2,\ldots n)$ as independent random variables uniformly distributed on $\T =\R/\Z$, given $n$ integer $\ge 2$ and $0<\delta<1$, the probability}
$$
P(M,n,\delta)= P(\forall m\in \{1,2,\ldots M\}\ \Big| \sum_{j=1}^n e^{2\pi i mx_j}\Big| <\delta n)
$$
\textit{satisfies}
$$
1-P (M,n,\delta) \le M \nu e^{-n\delta^2 a^2} \qquad (a^2 =\cos^2 \frac{\pi}{\nu})
$$
\textit{for all integers $\nu \ge 3$, and in particular for} $\nu =[\delta \pi \sqrt{2n}]$.

\vskip2mm

Therefore $P(M,n,\delta) >0$ when
$$
M < \frac{1}{\nu} e^{n\delta^2a^2} \qquad (a^2=\cos \frac{\pi}{\nu})\,.
$$
This is a better estimate than $\frac{1}{2} 2^{n\delta^{2}/4}$ already for $\nu=4$, and even for $\nu=3$ when $n\delta^2 \ge 8$.

\section{A question raised by Tur\'an and its discussion}

After stating and proving Theorem A (11.6) Tur\'an added the following comment:

``It would be interesting to find an explicit system $x_1,\ldots, x_n$ that satisfies the theorem''

Certainly Tur\'an had in mind a construction of the type of that he gave after Tijdeman of real $x_j$ $(j=1,2,\ldots n)$ such that
$$
\max_{1\le m \le n^A} \Big|\sum_{j=1}^n e^{2\pi im x_j}\Big| \le c(A)\sqrt{n\log n}\,.
$$
This construction involves prime numbers and gives an estimate of $c(A)$:
$$
c(A) = 6A+3\,, \qquad (\cite{Turan}\textrm{p}. 83)\,.
$$
It is not sufficient in order to obtain Theorem A but it provides a weaker result of the same type. Actually, taking
$$
\delta = c(A) \sqrt{\frac{\log n}{n}}\,,
$$
we have
$$
\begin{array}{ll}
A
	& \displaystyle = \frac{1}{6}\Big(\delta \sqrt{\frac{n}{\log n}} -3\Big)\,,   \\
\noalign{\vskip2mm}
n^A	&= \exp \Big(\dfrac{1}{6} (\delta \sqrt{n \log n} - 3 \log n)  \Big)\,.	
\end{array}
$$
Therefore, given $n\ge 3$ and $0<\delta<1$, we can exhibit $x_1,\ldots,x_n$ real such that (1) is valid for all integers $m$ such that
$$
1 \le m \le \exp \Big(\frac{1}{6} \delta \sqrt{n \log n} -3 \log n \Big)\,.
\leqno(3)
$$
Condition (3) is much stronger than (2), it means that the explicit construction in \cite{Turan}, pp.~82--83 gives a much weaker result than Theorem~A.

No such construction is available at the present time with (1) instead of (3). However, a random choice of the $x_j$ can be as good as an explicit construction from a practical point of view. It may provide good expamples with a high probability, since choosing
$$
M = \frac{1}{2} 2^{n \delta^{2/4}}
$$
in Theorem 1 gives
$$
M \nu e^{-n \delta^2 a^2} = \frac{\nu}{2} M^{-c}
$$
with
$$
c = \frac{4a^2}{\log 2} -1 \qquad \Big(a^2 = \cos^2 \frac{\pi}{\nu}\Big)\,.
$$
For instance $c>4,2$ when $\nu=10$: when $M$ is large we are nearly insured that Theorem~A works with a random choice of the $x_j$.

\section{A discrete version of Theorem 1}

In Theorem 1 the $x_j$ are independent random variables that are uniformly distributed on $\T$. Here we replace $\T$ by $\Z_N = \Z/N\Z$, $N$ prime $\ge 5$, and we choose $M=\frac{N-1}{2}$.

\vskip 2mm

\textsc{Theorem} 2.--- \textit{Considering the $X_j$ $(j=1,2,\ldots n)$ as independent random variables uniformly distributed on $\Z_N$, $N$ prime $\ge 5$, given $n$ integer $\ge 2$ and $0<\delta<1$, the probability}
$$
P(N,n,\delta) =P(\forall m \in \Z_N \ba \{0\}\ \Big| \sum_{j=1}^n e^{2\pi i m X_j/N} \Big| < \delta n\Big)
$$
\textit{satisfies}
$$
1 - P(N,n,\delta) < \frac{N-1}{2} \nu e^{-n \delta^2 a^2}\,, \qquad a^2 =\cos^2 \frac{\pi}{\nu}\,,
$$
\textit{for all integers $\nu \ge 3$ and in particular for $\nu = [\delta\pi\sqrt{2n}]$.}

\vskip2mm

The proof relies on the following lemma

\vskip2mm

\textsc{Lemma}.--- \textit{If $X$ is a random variable uniformly distributed on $\Z_N$, $N\ge 5$, then}
$$
E(e^{u\cos 2\pi X/N})\le e^{u^{2/4}} \qquad (u>0)\,.
\leqno(4)
$$

\vskip2mm

\textit{Proof of the lemma}. We compare the Taylor expansions on both sides:
$$
\begin{array}{l}
E(e^{u \cos 2\pi X/N}) = \displaystyle \sum_{k=0}^\infty \frac{u^k}{k!} J_k\\
\noalign{\vskip 2mm}
J_k=E\Big(\cos^k \dfrac{2\pi X}{N}\Big)\,.
\end{array}
$$
Using
$$
\cos^k u = 2^{-k} \sum_{j=0}^k \frac{k!}{j!(k-j)!} \cos (k-2j)u
$$
and
$$
E\Big(\cos \dfrac{2\pi jX}{N} \Big) =
\left\{\begin{array}{lll}
 1 &\mathrm{if} &j\in N\Z\\
0 &\mathrm{if} &j\notin N\Z
\end{array}\right. ,
$$
we see that
$$
\begin{array}{l}
J_0 =1,\ J_1=0,\ J_2=\dfrac{1}{2},\ J_3 =0,\ J_4=\dfrac{3}{8},\\
J_5=\dfrac{1}{16}\ \mathrm{or}\ 0\ \mathrm{and}\ J_6=\dfrac{1}{2}\  \mathrm{or}\ \dfrac{5}{6}\ \mathrm{according\ to\ } N=5 \ \mathrm{or}\ N\not= 5
\end{array}
$$
and
$$
J_k \le J_6 \le \frac{1}{2}\qquad \mathrm{if}\ k\ge 6\,.
$$
In both cases $u\le 1$ or $u>1$ we have
$$
\sum_{k=0}^\infty \frac{u^k}{k!} J_k \le 1+ \sum_{k=1}^\infty \frac{u^{2k}}{(2k)!}\Big(2k J_{2k-1} +J_{2k}+\frac{1}{2k+1}J_{2k+1}\Big)
$$
to be compared with
$$
e^{u^{2}/4} = 1 + \sum_{k=1}^\infty \frac{u^{2k}}{4^k k!}\,.
$$
It suffices to prove that for $k\ge 1$
$$
2k J_{2k-1} +J_{2k} + \frac{J_{2k+1}}{2k+1} \le \frac{(2k!)}{4^kk!} = \frac{(k+1)(k+2)\cdots (2k)}{4^k}\,.
$$
We check this inequality for $k=1,2,3$:
$$
\begin{array}{l}
0+\dfrac{1}{2} +0 \le \dfrac{2}{4}\\
\noalign{\vskip2mm}
0+\dfrac{3}{8} +\dfrac{1}{5\times 16} \le \dfrac{3\times 4}{4^2}\\
\noalign{\vskip2mm}
\dfrac{6}{16} +\dfrac{1}{2}+\dfrac{1}{7\times 2} \le \dfrac{4\times 5\times 6}{4^3}
\end{array}
$$
and for $k\ge 4$ we just use $J_\ell \le \frac{1}{2}$ $(\ell \ge 6)$ and we check
$$
\frac{1}{2}\Big(2k+1 +\frac{1}{2k+1}\Big) < \frac{(k+1)(k+2)\cdots (2k)}{4^k}\,.
$$
The lemma is proved.

Let us remark that (3) is not valid when $N=2$ or $N=3$.

The assumption that $N$ is prime is essential in the proof of the theorem. It insures that the $X_j$ and the $mX_j$ $(m=1,2,\ldots M)$ have the same joint distribution. The proof of Theorem 2 is copied from that of Theorem~1, with a slight modification: here again
$$
P\Big(\Big|\sum_{j=1}^n e^{2\pi i X_j/N}\Big| > \delta n\Big) \le \nu P \Big(\sum_{j=1}^n \cos \Big(\frac{2\pi X_j}{N} - \varphi\Big) > \delta na\Big)
$$
for some $\varphi$ multiple of $\frac{2\pi}{\nu}$ and $a=\cos \frac{\pi}{\nu}$, but now
$$
\begin{array}{l}
\displaystyle P\Big(\sum_{j=1}^n \cos \Big(\frac{2\pi X_j}{N}\! -\! \varphi\Big) > \delta na\Big) \le E \Big(\exp u \Big(\sum_{j=1}^n \cos \Big(\dfrac{2\pi X_j}{N} -\varphi\Big) -\delta na\Big)\Big)\hfill\\
\noalign{\vskip2mm}
= e^{-\delta nau}\Big(E \exp u\cos \Big(\dfrac{2\pi X}{N}-\varphi
\Big)\Big)^n
\end{array}
$$
and we shall check that this is
$$
\begin{array}{l}
\le e^{-\delta nau}\Big(E \exp u\cos \dfrac{2\pi X}{N}\Big)^n\\
\le e^{-\delta nau+nu^{2}/4}
\end{array}
$$
because of (4) and we proceed from that point as in Theorem~1. The main change is the inequality
$$
E \exp u\cos \Big(\frac{2\pi X}{N} -\varphi\Big) \le E \exp u \cos \frac{2\pi X}{N}
$$
and it is justified in expanding the first exponential in the form
$$
\sum_{j=0}^\infty \alpha_j(u) \cos j\Big(\frac{2\pi X}{N} - \varphi \Big)
$$

with $\alpha_j(u) \ge 0$ and using
$$
E \cos j\Big(\frac{2\pi X}{N} - \varphi\Big) \le E \cos j \frac{2\pi X}{N}
$$
for all $j$, multiples of $N$ or not.

Theorem 2 is proved and will be used in what follows.

\section{A new variation around the compressive sampling theorem of Cand\`es, Romberg and Tao}

Here we need new notations. We denote by $G$ the abelian group $\Z_N = \Z/N\Z$, considered as the space of times, $t$. The dual group $\hat{G}$, also represented  by $\Z_N$, is considered as the space of frequencies, $\omega$, and the duality is expressed~by
$$
(\omega,t) = e\Big(\frac{t\omega}{N}\Big)\,, \qquad e(u) = e^{2\pi iu}\,.
$$
A signal $x(t)$ is a complex--valued function on $G$. Its Fourier transform~is
$$
\hat{x}(\omega) = \frac{1}{\sqrt{N}} \sum_{t\in G} x(t) e\Big(\frac{-t\omega}{N}\Big)
$$
and $x(t)$ is reconstructed from the $\hat{x}(\omega)$ by the inversion formula
$$
x(t) = \frac{1}{\sqrt{N}} \sum_{\omega \in \hat{G}}
 \hat{x} (\omega) e\Big(\frac{t\omega}{N}\Big)\,.
 $$
 
The compressive sampling consists in reconstructing $x(t)$ by using only a small set of frequencies $\omega$, that is, $\omega\in \Omega$, $\Omega\subset \hat{G}$. It needs an assumption on the signal, usually that the signal is supported by a small set of times $t$, say $t\in S$, $S\subset G$. The cardinal of $S$ is denoted by $T : |S| =T$. Let us recall the notation:
$$
\|x\|_0 = \mathrm{cardinal\ of\ the\ support\ of\ }x\,.
$$
It is convenient to consider $x$ as an element of $\ell^1(G)$ and $\hat{x}$ as an element of $A(\hat{G})= \FF\ell^1 (G)$, with the norms
$$
\|\hat{x}\|_A = \|x\|_1 = \sum_{t\in G}|x(t)|\,.
$$

The process of reconstruction is as follows: $\hat{x}|_\Omega$ can be extended to $\hat{G}$ in many ways; if there is a unique extension of minimal norm in $A(\hat{G})$ and if this unique extension is $\hat{x}$, then $x$ is the solution of a problem of convex analysis (find a point of minimal norm in a given closed convex set) tractable from a numerical point of view. The reconstruction works under the assumption

\vskip2mm

\noindent$(\alpha)$:\hskip2mm $\hat{x}$ \textit{is the extension of $\hat{x}|_\Omega$ of minimal norm in} $A(\hat{G})$.

\vskip2mm

Here is the theorem of Cand\`es, Romberg and Tao (2006, \cite{CaRoTao,Cand}) that can be considered as a paradigm in the theory of compressive sampling, or compressed sensing.

\vskip2mm

\textsc{Theorem} B (2.1 in \cite{Cand}).--- \textit{Suppose that the signal $x$ is carried by a set of $T$ points $(\|x\|_0 \le T)$. Choose}
$$
f = [CT\, \log N] \qquad ([\ ]: \mathrm{integral\ point})
\leqno(5)
$$
\textit{and choose $\Omega$ randomly with the uniform distribution among all subsets of $\hat{G}$ such that $|\Omega|=f$.~If}
$$
C=22(1+\delta)
$$
\textit{the probability that $(\alpha)$ is valid satisfies}
$$
P((\alpha)) = 1 -O(N^{-\delta}) \qquad (N\rg \infty)
\leqno(6)
$$

I gave a number of variations about this theorem in \cite{JPK1} and \cite{JPK2}. Here is such a variation (just copied from the remark after (ii) in \cite{JPK2}):

\vskip2mm

\noindent $(\beta)$:\hskip2mm \textit{if there is a function}
$$
L(t) = \sum_{\omega\in \Omega} \hat{L}(\omega) e \big(\frac{t\omega}{N}\big)
$$
\textit{such that}
$$
\sup_{t\in G \ba \{0\}} |L(t)| < \frac{1}{2T} L(0)\,,
\leqno(7)
$$
\textit{then $(\alpha)$ is valid for all signals $x$ such that} $\|x\|_0 \le T$.

A more restricted form, used in \cite{JPK1} and \cite{JPK2}, is the following:

\vskip2mm

\noindent $(\gamma)$:\hskip2mm \textit{if\hskip4mm $K(t) = \sum\limits_{\omega\in \Omega} e\big(\frac{t\omega}{N}\big)$ satisfies}
$$
\sup_{t\in G \ba\{0\}} |K(t)| < \frac{1}{2T} K(0)\,,
\leqno(8)
$$
\textit{then $(\alpha)$ is valid for all signals $x$ such that} $\|x\|_0 \le T$.

\vskip2mm

How to construct such a function $L(t)$, or such a function $K(t)$? Here again a random choice is convenient, and it leads to the following result (cf. $V2''$ in~\cite{JPK2}):

\vskip2mm

\textsc{Theorem} C.--- \textit{Suppose that $\Omega$ is produced by a random selection on $\hat{G}$ with parameter $\tau$ $(0<\tau<1)$, meaning that the events $(\omega\in \Omega)$ are independent and $P(\omega\in \Omega)=\tau$ for all $\omega\in \hat{G}$ (then $|\Omega|$ has a binomial distribution $B(\tau,N)$ with mean value $\tau N$ and variance $\tau(1-\tau)N)$. Assume that}
$$
\tau N = 4C(T^2+1) \log N\qquad (C>1)
\leqno(9)
$$
\textit{and $N=\pm1$ modulo $6$. Then}
$$
P((8)) > 1 - \nu N^{-Ca^2+1}
$$
\textit{for all integers $\nu\ge 3$ and $a^2=\cos^2\frac{\pi}{\nu}$. Consequently}
$$
P(\forall x : \|x\|_0 \le T \ \ (\alpha)) = 1-O(N^{-\delta}) \qquad (N\rg \infty)
\leqno(10)
$$
\textit{for all} $\delta <C-1$.

\vskip2mm

The assumption (9) is much stronger than (5), but the result (10) is much stronger than~(6).

Assuming that $N$ is prime, Theorem 2 is a new way to look at Theorem~C. Now the $X_j$ of Theorem~2 are random frequencies, not necessarily distinct, and $\Omega$ is the random set of values taken by the~$X_j$, $j=1,2,\ldots n$, therefore $|\Omega|\le n$. We can write
$$
\sum_{j=1}^n e^{2\pi im X_j/N} = \sum_{\omega\in \Omega} \hat{L}(\omega)e\Big(\frac{m\omega}{N}\Big) = L(m)
\leqno(11)
$$
and, choosing $\delta=\frac{1}{2T}$, $P((7))$ is the same as $P(N,n,\delta)$ in Theorem~2. The result reads as follows

\vskip2mm

\textsc{Theorem} 3.--- \textit{Suppose that $\Omega$ is the set of values taken by $n$ independent random variables $X_j$ uniformly distributed on $\Z_N$ and assume that}
$$
n = 4CT^2 \log N \qquad (C>1)
\leqno(12)
$$
\textit{and $N$ is prime $\ge 5$. Then the function $L$ defined in (11) satisfies}
$$
P((7)) > 1 - \frac{\nu}{2} N^{1-Ca^2}
\leqno(13)
$$
\textit{for all integers $\nu\ge 3$ and $a^2=\cos^2 \frac{\pi}{\nu}$. Consequently}
$$
P(\forall x : \|x\|_0 \le T\ \ (\alpha)) = 1 - O(N^{-\delta}) \qquad (N\rg \infty)
\leqno(14)
$$
\textit{for all} $\delta<C-1$.

\vskip2mm

The conclusion (14) is the same as (10) and the assumption (12) is only a slight improvement of (9). However, using (13) in an optimal way, it is interesting to apply Theorem~3 when $T$, small, and $N$ are given. Forgetting the fact that $\nu$ is an integer, the best choice of $\nu$ in (13) satisfies
$$
C \frac{\pi}{\nu} \sin \frac{2\pi}{\nu} \log N=1\,.
$$
We choose a slight variation:
$$
\begin{array}{c}
\sin^2\dfrac{\pi}{\nu} = \dfrac{1}{2C\log N}\,, \ \ a^2 =1-\dfrac{1}{2C\log N}\,,\\
\noalign{\vskip2mm}
\nu < \pi \sqrt{2C\log N}
\end{array}
$$
\textit{That gives}
$$
P((7)) > 1 - \frac{\pi}{2} \sqrt{2Ce\log N} N^{1-C}
\leqno(15)
$$
\textit{instead of $(13)$ and accordingly}
$$
P(\forall x : \|x\|_0 \le T (\alpha)) > 1 - \frac{\pi}{2} \sqrt{2Ce\log N} N^{1-C}
\leqno(16)
$$
\textit{much more precise as} (14). As an example, let me choose $N=997$ and $T=2$, and estimate $n$ and $p=\frac{\pi}{2}\sqrt{2Ce\log N} N^{1-C}$ for a few values of~$C$:
$$
\begin{array}{lll}
C=2      &\qquad n=242    	&\qquad p<0,036\\
C=3      &\qquad n=332		&\qquad p<0,000044
\end{array}
$$
The random choice of $\Omega$ with $n=332$ is near to insure that $(\alpha)$ is valid for  all signals carried by 2 points in $\Z_{997}$. Though that has no practical value it shows the power of random methods compared to explicit constructions (here not yet discovered).

Let me conclude with an observation on random subsets of a given finite set, here $\Z_N$. In Theorems B, C, 2 and 3 I used $P(\ )$ for the probability of an event consisting in a collection of sets $\Omega$ included in $\Z_N$, but the probability spaces were different. Here are the probability spaces we met.

1 (Theorem B) $P(\ )$ is the uniform distribution on all subsets consisting of $f$ points.

2 (Theorem C) $P(\ )$ is defined by the condition that all events $(\omega\in \Omega)$ are independent and have the same probability $\tau$ (random selection).

3 (Theorems 2 and 3) $P(\Omega)$ is the probability that $\Omega$ is exactly the range of $X_1,X_2,\ldots X_n$, independent random variables uniformly distributed on~$\Z_N$.

Another natural probability space is

4 (Poisson point process) $P(\Omega)$ is a Poisson variable of parameter~$\tau|\Omega|$.

Roughly speaking, $f$, $\tau N$ and $n$ have the same role, but the relation between these probabilities deserves attention. All of them are invariant under any permutation of the given set (here $\Z_N$). In other words, the conditional probability when $|\Omega|$ is given is uniform. They are well defined by the distribution of $|\Omega|$: the Dirac measure at $f$ in case 1, the Bernoulli distribution $B(\tau,N)$ in case~2, the Poisson distribution $P(\tau)$ in case~4, and interesting distribution carried by $\{1,2,\ldots n\}$ in case~3. The comparison between  these distributions of $|\Omega|$ is done in \cite{CaRoTao} and \cite{JPK2} for cases~1 and 2, it is easy to extend it to case~4, the case~3 needs more attention. It is known among probabilists as the occupation problem and much is known about (see for example \cite{Fel} or \cite{Dur}); I thank Gregory Miermont for this observation.

\eject


\begin{thebibliography}{3}

\bibitem{Cand}
J. \textsc{Cand\`es}.--- \textit{Compressive sampling,}  Proceedings of the International Congress of Mathematicians, Madrid 2006.

\bibitem{CaRoTao}
J. \textsc{Cand\`es}, J. \textsc{Romberg} and T. \textsc{Tao}.--- \textit{Robust Uncertainty Principles : Exact Signal Reconstruction from Highly Incomplete Frequency Information}, IEEE Transactions on Information Theory 20, 2 (2006), 489--509.

\bibitem{Dur}
R. \textsc{Durrett}.--- \textit {Probability: Theory and Examples,}  Wedsworth and Brooks/Cole, Pacific Grove, Calif, 1991.

\bibitem{Fel}
W. \textsc{Feller}.--- \textit{An Introduction to Probability Theory and Its Applications,}  Volume~I, Third edition, Wiley, New York 1967.

\bibitem{JPK1}
J.--P. \textsc{Kahane}.--- \textit{Idempotents et \'echantillonnage parcimonieux,}  C.R. Acad. Sci.

\bibitem{JPK2}
J.--P. \textsc{Kahane}.--- \textit{Variantes sur un th\'eor\`eme de Cand\`es, Romberg et Tao,}  Ann. Inst. Fourier.

\bibitem{Turan}
P. \textsc{Tur\'an}.--- \textit{On a New Method of Analysis and its Applications,}  Pure and Applied Mathematics, John Wiley 1984.


\vskip4mm

\begin{tabular}{p{6cm}l}
&Jean--Pierre Kahane,\\
&Laboratoire de Math\'ematiques\\
&Universit\'e Paris--Sud \`a Orsay\\
&Jean-Pierre.Kahane@math.u-psud.fr
\end{tabular}





\end{thebibliography}
\end{document}